\catcode`\à=\active \def à{\` a}
\catcode`\â=\active \def â{\^ a}
\catcode`\é=\active \def é{\' e}
\catcode`\è=\active \def è{\` e}
\catcode`\ê=\active \def ê{\^ e}
\catcode`\ë=\active \def ë{\" e}
\catcode`\ù=\active \def ù{\` u}
\catcode`\û=\active \def û{\^ u}
\catcode`\ô=\active \def ô{\^ o}
\catcode`\ç=\active \def ç{\c c}
\catcode`\î=\active \def î{\^ \i}
\catcode`\ï=\active \def ï{\" \i}

\magnification=1250

\input color

\vsize=190mm

\baselineskip 11.5pt

\font\pg=cmss10 scaled 900

\font\sc=cmcsc10 at 9pt

\font\zazou = cmb10 at 12pt

\font\zazoue = cmb10 at 7pt

\font\pg=cmss10 scaled 900

\font\tfp=cmssbx10 scaled 1200

\font\tfpp=cmssbx10 scaled 1000




\font\mathb=msbm10 scaled 1000

\let\Bbb=\mathbbb

\input amssym

\input epsf

\font\ib = cmbxsl10

\centerline{\zazoue \'EQUATION COHOMOLOGIQUE DU FLOT AFFINE DE}

\centerline{\zazoue REEB SUR LA VARI\'ET\'E DE HOPF ${\Bbb S}^n\times {\Bbb S}^1$}

\medskip
\centerline{{\sc Aziz El Kacimi Alaoui}}
\medskip
\centerline{(Septembre 2019)}
 \vskip0.5cm
\noindent {\parindent=0.8cm\narrower {\bf Résumé.} {\pg On calcule explicitement la cohomologie feuilletée $H_{\cal F}^\ast (M)$
du flot affine de Reeb ${\cal F}$ sur la variété de Hopf ${\Bbb S}^n\times {\Bbb S}^1$. L'espace vectoriel $H_{\cal F}^1(M)$ contient exactement les
obstructions à la résolution de l'équation cohomologique $X\cdot f=g$ et permet de décrire l'espace des distributions invariantes par tout champ $X$
non singulier définissant le feuilletage ${\cal F}$.
\par }}
\vskip0.5cm

{\zazou 1. Rappels}
\medskip
\noindent  Dans cette section, on définit les principales notions considérées dans ce texte : celle de feuilletage, de cohomologie
feuilletée, de suite spectrale d'un revêtement...Tous les objets géométriques (variétés, fonctions, champs de vecteurs, formes différentielles...) seront supposés
de classe $C^\infty $. Et sauf mention expresse du contraire, les fonctions et les formes différentielles seront à valeurs complexes.
\smallskip
\noindent {\tfpp 1.1. D\'efinion.} {\it Soit $M$ une vari\'et\'e (connexe) de dimension $m+n$. Un
{\ib feuilletage} ${\cal F}$ de {\ib codimension} $n$ (ou de {\ib dimension} $m$) sur $M$
est la donn\'ee d'un recouvrement ouvert ${\cal U}=\{ U_i\}_{i\in I}$ et,
pour tout $i$, d'un diff\'eomorphisme $\varphi_i: {\Bbb R}^{m+n}\longrightarrow U_i$ tel que, sur toute intersection non vide
$U_i\cap U_j$, le diff\'eomorphisme de changement de coordonn\'ees $\varphi_j^{-1}\circ \varphi_i: (x,y)\in \varphi_i^{-1}(U_i\cap U_j)
\longrightarrow(x',y')\in \varphi_j^{-1}(U_i\cap U_j)$
soit de la forme $x'=\varphi_{ij}(x,y)$ et $y'=\gamma_{ij}(y)$.}

\smallskip

On peut aussi voir ${\cal F}$ comme la donn\'ee d'un sous-fibr\'e $\tau $ de
rang $m$ du fibr\'e tangent $TM$ {\it compl\`etement int\'egrable},
c'est-\`a-dire, pour toutes sections $X,Y\in C^\infty (\tau )$ de
$\tau $ ({\it i.e.} des champs de vecteurs sur $M$ tangents \`a
$\tau $), le crochet $[X,Y]$ est encore une section de $\tau $. Les
sous-vari\'et\'es connexes tangentes \`a $\tau $ sont appel\'ees
{\it feuilles} de ${\cal F}$.

\smallskip

\noindent {\tfpp 1.2. Exemples}
\smallskip

Ce n'est pas ce qui manque, la théorie des feuilletages en est riche. Mais nous n'en donnerons que
deux, à la fois assez importants et de description relativement facile.

\smallskip
\noindent {\tfpp i) Le feuilletage lin\'eaire sur le tore ${\Bbb T}^2$.} Soit  $\widetilde M={\Bbb R}^2$ et consid\'erons l'\'equation diff\'erentielle lin\'eaire
$dy-\alpha dx=0$ o\`u  $\alpha $ est un nombre r\'eel.
Elle a pour solution g\'en\'erale
$y=\alpha x+c$ avec $c\in {\Bbb  R}$ ; cette solution est une fonction lin\'eaire
dont le graphe est une droite $F_c$. Lorsque $c$
varie dans ${\Bbb R}$, on obtient une famille de droites parall\`eles qui remplit le plan $\widetilde M$
et qui y  d\'efinit un feuilletage   $\widetilde{{\cal F}}$ (ses feuilles sont les droites $F_c$). L'action
 naturelle de ${\Bbb Z}^2$ sur $\widetilde M$ pr\'eserve
$\widetilde{{\cal F}}$ ({\it i.e.}  l'image de chaque feuille
$\widetilde{{\cal F}}$ par une translation enti\`ere est une feuille de
$\widetilde{{\cal F}}$). Le feuilletage  $\widetilde{{\cal F}}$ induit alors un feuilletage
${\cal F}$ sur le tore ${\Bbb T}^2={\Bbb R}^2/{\Bbb
Z}^2$. Les feuilles de ${\cal F}$ sont toutes diff\'eomorphes au cercle  ${\Bbb S}^1$
si $\alpha $ est rationnel et \`a la droite r\'eelle sinon  (comme on le voit sur le dessin ci-dessous). En fait si
$\alpha $ n'est pas rationnel,
toute feuille de  ${\cal F}$ est dense ; ceci montre que m\^eme si, localement, un feuilletage
est simple, sa structure globale peut \^etre un peu compliqu\'ee.

\smallskip
\centerline{\epsfxsize=6cm \epsfbox{Tore.eps}}
\smallskip

\newif\ifpremierepage\premierepagefalse
\def\makeheadline{%
\ifpremierepage\global\premierepagefalse \else \vbox to
0pt{\vskip-22.5pt \line{\vbox to 8.5pt{}\ifodd\pageno \textedroite
\else \textegauche \fi } \vskip4pt
}\nointerlineskip \fi}
\def\boxit#1#2{\setbox1=\hbox{\kern#1{#2}\kern#1}%
\dimen1=\ht1 \advance\dimen1 by #1 \dimen2=\dp1 \advance\dimen2 by
#1
\setbox1=\hbox{\vrule height\dimen1 depth\dimen2\box1\vrule width 1mm}%
\setbox1=\vbox{\hrule\box1\hrule height 0.5mm depth 0.5mm}%
\advance\dimen1 by .4pt \ht1=\dimen1 \advance\dimen2 by .4pt
\dp1=\dimen2 \box1\relax}

\gdef\textegauche{{ } \hfil {\pg  A. El Kacimi}\hfil }
\gdef\textedroite{{ }\hfil {\pg Cohomologie feuillet\'ee du flot affine de Reeb}\hfil }

\noindent {\tfpp ii) Le feuilletage de Reeb sur la $3$-sph\`ere.} Soit $M$ la   sph\`ere ${\Bbb S}^3$ vue comme suit dans ${\Bbb C}^2$ : ${\Bbb S}^3=\{ (z_1,z_2)\in
{\Bbb C}^2:\vert z_1\vert^2+\vert z_2\vert^2=1\} $. On note ${\Bbb
D}$ le disque unit\'e ouvert de ${\Bbb C}$ et $\overline{\Bbb D}$ son adh\'erence (qui est le disque unit\'e ferm\'e $\{ z\in {\Bbb C}:\vert z\vert
\leq 1\} $). Les deux ensembles $M_+=\left\{ (z_1,z_2)\in {\Bbb S}^3:\vert z_1\vert^2\leq
{{1}\over {2}}\right\} $ et $M_-=\left\{ (z_1,z_2)\in
{\Bbb S}^3:\vert z_2\vert^2\leq  {{1}\over {2}}\right\}
$ sont diff\'eomorphes \`a  $\overline{\Bbb D}\times {\Bbb S}^1$, ont
le tore  ${\Bbb T}^2$ comme bord commun :
$${\Bbb T}^2=\partial M_+=\partial M_-=
 \left\{ (z_1,z_2)\in {\Bbb S}^3:\vert z_1\vert^2 =\vert z_2\vert^2=
{{1}\over {2}}\right\} $$ et leur  r\'eunion est  ${\Bbb S}^3$.
Il est alors clair que  ${\Bbb S}^3$ peut \^etre obtenue en recollant $M_+$ et
$M_-$ le long de leur bord commun  par le  diff\'eomorphisme
 $(z_1,z_2)\in \partial
M_+\longrightarrow (z_2,z_1)\in \partial M_-$ {\it i.e.}  on identifie
$(z_1,z_2)$ \`a   $(z_2,z_1)$ dans la r\'eunion disjointe  $M_+\coprod
M_-$ : on recolle deux tores solides en identifiant un m\'eridien du bord du premier
\`a un parall\`ele du bord du second. Soit $f:{\Bbb D} \longrightarrow {\Bbb R}$ la  fonction
d\'efinie par : $$f(z)=\hbox{exp}\left( {1}\over {1-\vert
z\vert^2}\right) .$$ Notons   $t$ la deuxi\`eme coordonn\'ee dans ${\Bbb
D}\times {\Bbb R}$. La famille de surfaces $(S_t)_{t\in {\Bbb R}}$
obtenue en translatant le graphe  $S$ de $f$  le long de l'axe des $t$
d\'efinit un feuilletage sur ${\Bbb D}\times {\Bbb R}$. Si on rajoute le cylindre   ${\Bbb S}^1\times {\Bbb R}$, o\`u ${\Bbb S}^1$ est vu
comme le bord de $\overline{\Bbb D}$, on obtient un feuilletage   $\widetilde{\cal  F}$ de codimension $1$ sur  $\overline{\Bbb D}\times {\Bbb
R}$.

\medskip

\centerline{\hskip0.4cm\epsfxsize=14cm \epsfbox{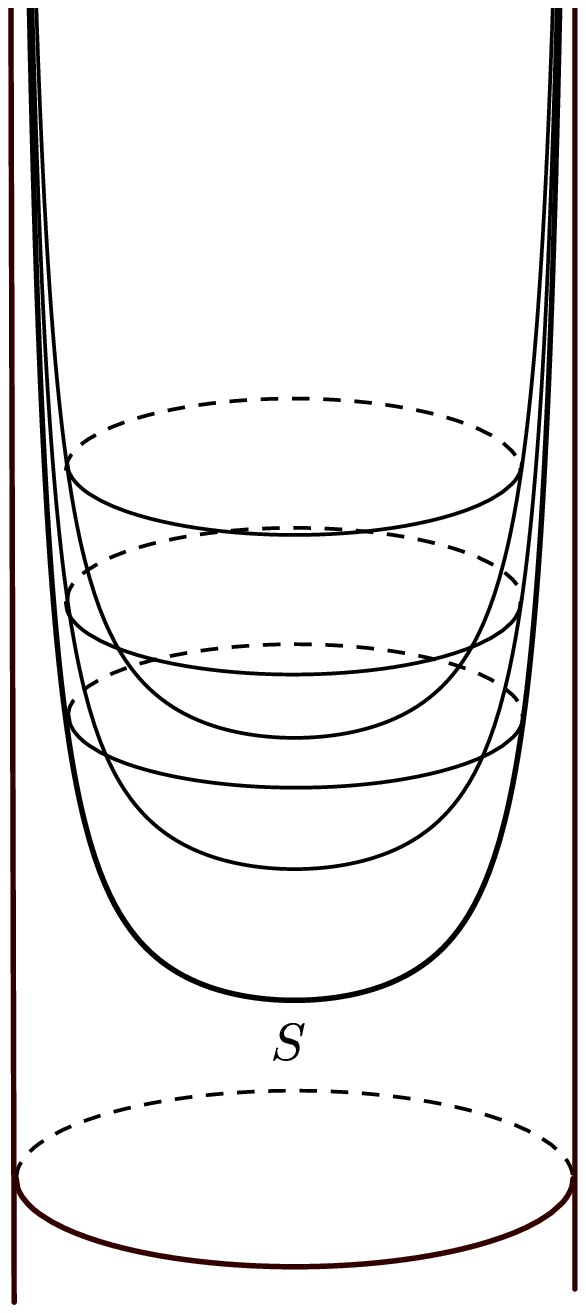}}

Par construction même, le feuilletage $\widetilde{{\cal F}}$ est invariant par les
transformations $(z,t)\in \overline{\Bbb D}\times {\Bbb R}\longmapsto
(z,t+1)\in \overline{\Bbb D}\times {\Bbb R}$ et induit donc un feuilletage
  ${\cal F}_0$ sur le l'espace quotient :
$$\overline{\Bbb D}\times {\Bbb R}/(z,t)\sim (z,t+1)\simeq
\overline{\Bbb D}\times {\Bbb S}^1.$$
Il a le bord ${\Bbb T}^2={\Bbb S}^1\times {\Bbb S}^1$
comme feuille compacte. Toutes les autres sont diff\'eomorphes \`a    ${\Bbb R}^2$ ({\it cf.} dessin ci-dessous).

\centerline{\epsfxsize=4cm \epsfbox{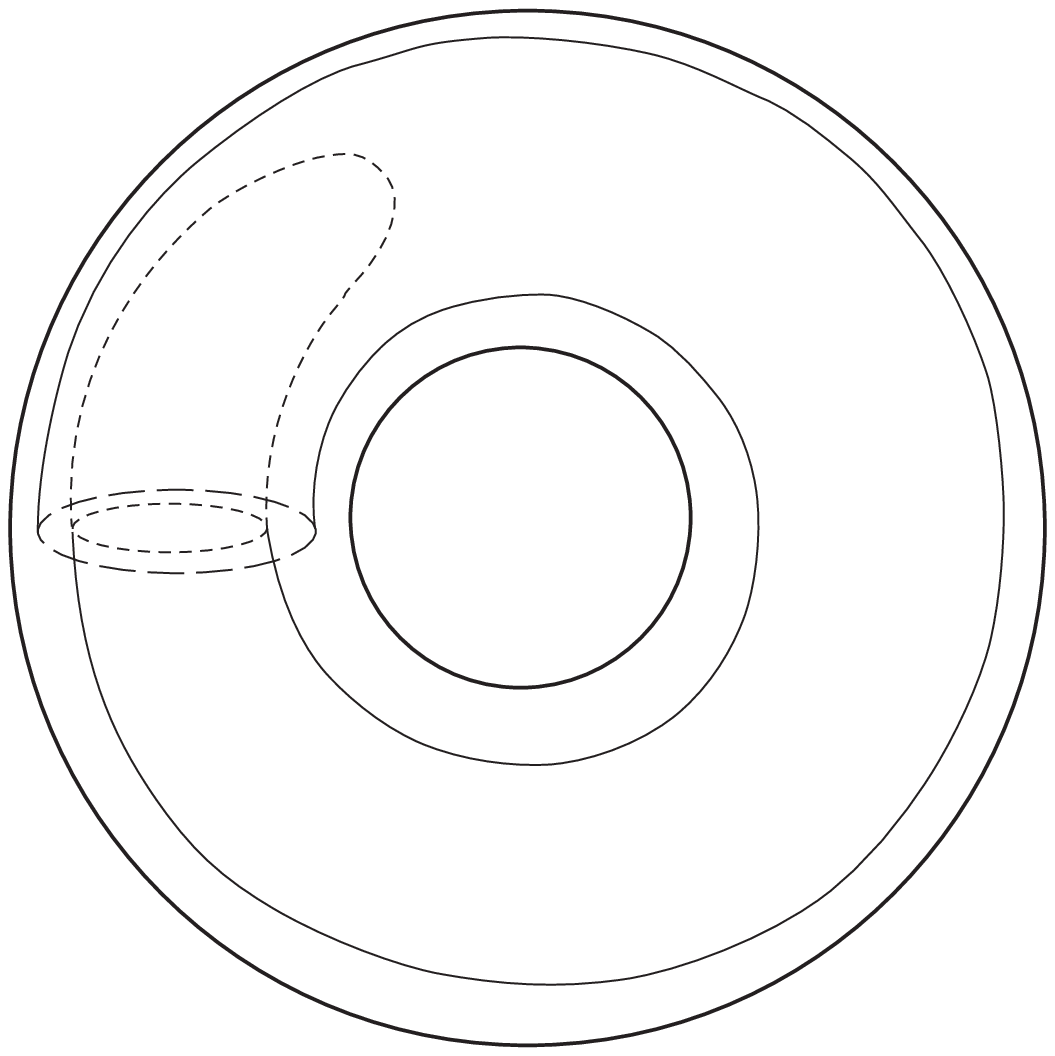}}
\medskip
\noindent Comme  les morceaux $M_+$ et $M_-$ sont diff\'eomorphes \`a
$\overline{\Bbb D}\times {\Bbb S}^1$, ${\cal F}_0$ d\'efinit sur $M_+$
et $M_-$  respectivement  deux feuilletages ${\cal F}_+$ et ${\cal
F}_-$ ; leur  recollement le long de la feuille compacte ${\Bbb T}^2$ donne un feuilletage ${\cal F}$ sur ${\Bbb
S}^3$ appel\'e {\it feuilletage de Reeb}.  Toutes les feuilles sont diff\'eomorphes au plan   ${\Bbb R}^2$
\`a l'exception de celle qui provient des bords de   $M_+$ et
$M_-$ qui est diff\'eomorphe au tore ${\Bbb T}^2$.

\smallskip

Le choix de la fonction $f(z)=\hbox{exp}\left( {1}\over {1-\vert
z\vert^2}\right) $ pour définir le feuilletage intermédiaire sur ${\Bbb D}\times {\Bbb R}$ par les surfaces $S_t$ n'est pas anodin. Le fait
que son graphe soit infiniment tangent à la feuille ${\Bbb S}^1\times {\Bbb R}$ du bord de $\overline{\Bbb D}\times {\Bbb R}$ permet de recoller
les feuilletages ${\cal F}_+$ et ${\cal F}_-$ de manière $C^\infty $ en le feuilletage de Reeb sur ${\Bbb S}^3$ (mais qui ne saurait être analytique
par le théorème de Haefliger).

\smallskip

La référence [Go] est une bonne présentation  de tout ce qui
a été fait en théorie des feuilletages depuis sa naissance jusque à peu près 1990,
avec une bibliographie bien fournie et une préface de Georges Reeb.

\medskip
\noindent {\tfpp 1.3. Cohomologie feuilletée}
\medskip
$\bullet $ Soit ${\cal F}$ un feuilletage de dimension $m$ sur $M$. On note $T{\cal F}$ son fibré tangent et $\nu $ un
sous-fibré supplémentaire dans $TM$, par exemple le fibré orthogonal pour une métrique riemannienne.
On dira qu'une forme différentielle $\alpha $ de degré $\ell $ sur $M$ est {\it feuilletée} si $\alpha (X_1,\cdots ,X_\ell )=0$
lorsque l'un des vecteurs $X_1,\cdots ,X_\ell $ est tangent à $\nu $. Soit $\Omega_{\cal F}^\ell (M)$ l'espace vectoriel de telles formes.
Pour $\ell =0$, $\Omega_{\cal F}^0(M)$ est l'espace $C^\infty (M)$ des fonctions $C^\infty $ sur $M$.   On a un
op\'erateur de diff\'erentiation ext\'erieure le long des feuilles : $\Omega_{\cal F}^\ell (M)\buildrel {d_{\cal F}}\over \longrightarrow
\Omega_{\cal F}^{\ell +1}(M)$ d\'efini (comme dans le cas classique), pour $X_1,\cdots ,X_\ell $ tangents à ${\cal F}$, par
la formule :
$$\eqalign{d_{\cal F}\alpha (X_1,\cdots ,X_{\ell +1})&=\sum_{i=1}^{\ell +1}(-1)^iX_i\cdot
\alpha (X_1,\cdots ,\widehat X_i,\cdots ,X_{\ell +1})\cr
&+\sum_{i<j}(-1)^{i+j}\alpha ([X_i,X_j],X_1,\cdots ,\widehat
X_i,\cdots ,\widehat X_j,\cdots ,X_{\ell +1})}$$ o\`u $\widehat X_i$
signifie qu'on a omis l'argument $X_i$. On v\'erifie facilement que
$d_{\cal F}$ est de carr\'e nul. On obtient ainsi un
complexe diff\'erentiel (dit {\it complexe feuillet\'e} de $(M,{\cal F})$) :
$$0\longrightarrow \Omega_{\cal F}^0(M)\buildrel {d_{\cal F}}\over
\longrightarrow \Omega_{\cal F}^1(M)\buildrel {d_{\cal F}}\over
\longrightarrow \cdots  \buildrel {d_{\cal F}}\over \longrightarrow
\Omega_{\cal F}^{m-1}(M)\buildrel {d_{\cal F}}\over \longrightarrow
\Omega_{\cal F}^m(M) \longrightarrow 0. \leqno{(1)}$$
Sa cohomologie (en tout degré $\ell $) :
$$H_{\cal F}^\ell (M)=Z_{\cal F}^\ell (M)/B_{\cal F}^\ell (M)\leqno{(2)}$$
ne dépend pas du choix du fibré supplémentaire $\nu $ ; on l'appelle
$\ell ^{\hbox{\`eme}}$ espace vectoriel de {\it cohomologie
feuillet\'ee} de $(M,{\cal F})$.  C'est un invariant important du
feuilletage. Par exemple, le dual topologique de l'espace $H_{\cal F}^m(M)$
contient les {\it cycles feuillet\'es} au sens de [Su] et donc,  en
particulier,  les mesures transverses invariantes ({\it cf.} [Ek]).
Le calcul de $H_{\cal F}^\ast (M)$ est souvent tr\`es ardu.   Pour
une utilisation int\'eressante de la cohomologie feuillet\'ee dans
l'\'etude de la rigidit\'e de certaines actions de groupes de Lie
voir [MM].

\smallskip
 $\bullet $ Une fonction $f$ sur $M$ est dite ${\cal F}$-{\it basique} (ou simplement {\it basique}) si elle  v\'erifie $d_{\cal F}f=0$, c'est-à-dire si elle est constante sur les feuilles.
Le faisceau des germes de telles fonctions sera not\'e ${\cal O}_{\cal F}$ ; il admet
une r\'esolution fine ({\it cf.} [Va]) :
$$0\longrightarrow {\cal O}_{\cal F}\buildrel {d_{\cal F}}\over
\longrightarrow {\Omega}_{\cal F}^0\buildrel {d_{\cal F}}\over
\longrightarrow {\Omega}_{\cal F}^1\buildrel {d_{\cal F}}\over
\longrightarrow \cdots  \buildrel {d_{\cal F}}\over \longrightarrow
{\Omega}_{\cal F}^{m-1}\buildrel {d_{\cal F}}\over \longrightarrow
{\Omega}_{\cal F}^m \longrightarrow 0 \leqno{(3)}$$
o\`u ${\Omega}_{\cal F}^\ell $ est le faisceau des germes des $\ell $-formes feuillet\'ees. Comme
$\Omega_{\cal F}^\ell (M)$ est l'espace des sections globales de ${\Omega}_{\cal F}^\ell $, on a un isomorphisme canonique :
$$H_{\cal F}^\ell (M)\simeq H^\ell (M,{\cal O}_{\cal F}).\leqno{(4)}$$
Les deux d\'efinitions permettent de faire des calculs, chacune a son utilité suivant la nature des exemples et la mani\`ere dont ils sont d\'ecrits.
Mais c'est plutôt la version ``faisceau" que nous utiliserons  le plus  dans la suite.

\smallskip

\noindent {\tfpp 1.4. Suite spectrale d'un ${\Bbb Z}$-revêtement feuilleté}
\smallskip
Soit $\widetilde M$ une variété munie d'un feuilletage $\widetilde{\cal F}$. Un {\it automorphisme} de $\widetilde{\cal F}$ est un
difféomorphisme $\widetilde M\longrightarrow \widetilde M$ envoyant feuille sur feuille. L'ensemble Aut$(\widetilde{\cal F})$ de
tels automorphismes est un groupe pour la composition des applications.

Une {\it action fidèle} du groupe ${\Bbb Z}$ sur $(\widetilde M,\widetilde{\cal F})$ est la donnée d'un morphisme injectif de
groupes ${\Bbb Z}\longrightarrow \hbox{Aut}(\widetilde{\cal F})$,
donc de l'image $\gamma $ de l'élément $1$ (générateur du groupe ${\Bbb Z}$). Si cette action est libre et propre, le quotient $M=\widetilde M/\langle \gamma \rangle $
est une variété de même dimension que $\widetilde M$ ;  $\widetilde{\cal F}$ y induit alors un feuilletage ${\cal F}$ de même dimension (et donc de même codimension).
La projection canonique $\pi :\widetilde M\longrightarrow M$ est un revêtement feuilleté {\it i.e.} un revêtement qui, en plus,
envoie toute feuille de $\widetilde{\cal F}$ sur une feuille de ${\cal F}$.

De ce qui précède on déduit que l'image réciproque $\pi^\ast \left( {\cal O}_{{\cal F}}\right) $ par $\pi $ du faisceau
${\cal O}_{{\cal F}}$ est le faisceau ${\cal O}_{\widetilde{\cal F}}$. On a donc une suite spectrale de terme :
$$E_2^{pq}=H^p(\Gamma ,H^q(\widetilde M,{\cal O}_{\widetilde{\cal F}}))\leqno{(5)}$$
et convergeant vers $H^{p+q}(M,{\cal O}_{\cal F})$. (L'espace vectoriel $H^q(\widetilde M,{\cal O}_{\widetilde{\cal F}})$ est vu comme $\Gamma $-module
pour l'action induite par celle de $\Gamma =\langle \gamma \rangle $ sur $\widetilde M$.) Comme $\Gamma
\simeq {\Bbb Z}$, la diff\'erentielle $d_2:E_2^{pq}\longrightarrow
E_2^{p+2,q-1}$ est nulle et donc la suite spectrale stationne au
terme $E_2$. Ce qui donne, pour tout $\ell \in {\Bbb N}$, $H_{\cal F}^\ell (M)=\displaystyle E_2^{0\ell }\oplus E_2^{1,\ell-1}.$
En particulier :
$$H_{\cal F}^0(M)=H^0(\Gamma, H^0(\widetilde{M},{\cal O}_{\widetilde{\cal F}}))=\{ \hbox{fonctions basiques $\Gamma $-invariantes sur $\widetilde M$}\} .\leqno{(6)}$$
$$H_{\cal F}^1(M)=E_2^{01}\oplus E_2^{10}=H^0(\Gamma, H^1(\widetilde{M},{\cal O}_{\widetilde{\cal F}}))
\oplus H^1(\Gamma, H^0(\widetilde{M},{\cal O}_{\widetilde{\cal F}})).\leqno{(7)}$$
Si en plus $\widetilde{\cal F}$ est un flot (et donc aussi ${\cal F}$), les espaces
vectoriels $H^\ast (\widetilde{M},{\cal O}_{\widetilde{\cal F}})$ et $H^\ast (M,{\cal O}_{\cal F})$ sont nuls pour $\ast \geq 2$ et donc la cohomologie
feuilletée $H^\ast (M,{\cal O}_{\cal F})$ se réduit à celle donnée par les formules (6) et (7). Ce sera le cas pour le flot affine de Reeb, objet de cette note.
\medskip
{\zazou 2. Construction du flot affine de Reeb}
\medskip
\noindent  Pour $n\geq 2$  entier naturel, on note ${\Bbb E}$ l'espace euclidien ${\Bbb R}^n$ et $(z,t)$  les coordonnées d'un point quelconque de ${\Bbb E}\times {\Bbb R}$. On pose
${\Bbb E}^\ast ={\Bbb E}\setminus \{0\} $ et $\widetilde M={\Bbb E}\times {\Bbb R}\setminus \{ (0,0)\} $.
\smallskip
\noindent {\bf 2.1.}  Le système différentiel $dz=0$ définit un feuilletage (un flot)
$\widetilde{\cal F}$ sur $\widetilde M$ ; ses feuilles sont les courbes intégrales du champ de vecteur $\widetilde T={\partial \over {\partial t}}$.
De la même manière, l'équation différentielle $dt=0$ définit un feuilletage $\widetilde{\cal V}$ de codimension $1$.
\smallskip
-- Toute feuille de $\widetilde{\cal F}$ passant par $(z,t)$ avec $z\neq 0$ est une droite ; les autres sont deux demi-droites $\widetilde F_+=\{ (0,t):t>0\} $ et
$\widetilde F_-=\{ (0,t):t<0\} $.

-- Toute feuille de $\widetilde{\cal V}$ passant par $(z,t)$ avec $t\neq 0$ est une copie de l'espace ${\Bbb E}$. Celle qui reste (qui passe par les $(z,0)$), qu'on notera
$\widetilde V$, est isomorphe à
${\Bbb E}^\ast $.
\medskip
\centerline{\hskip0.4cm\epsfxsize=10.5cm
\epsfbox{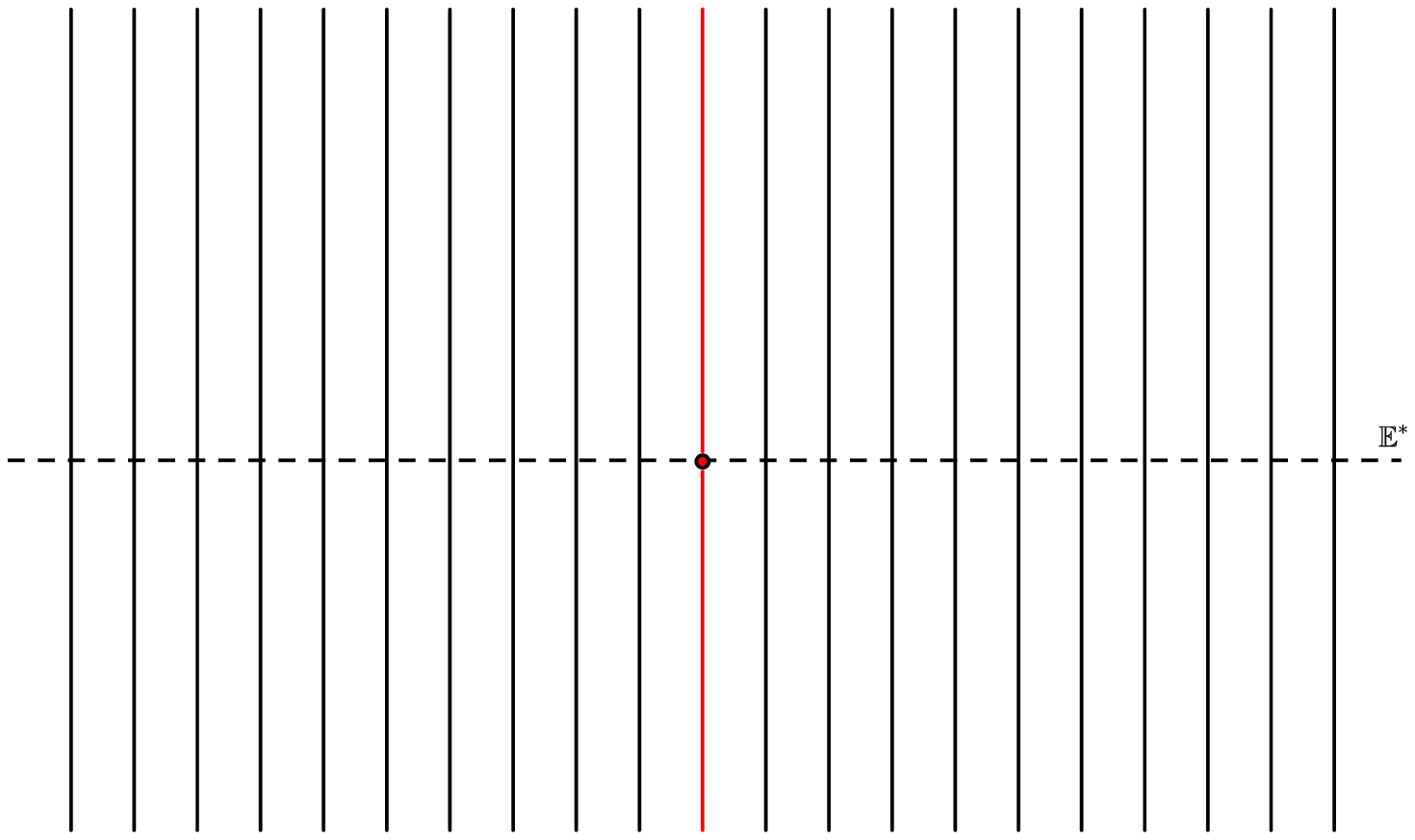}}
\smallskip
\noindent {\bf 2.2.} Soient $\lambda \in ]0,1[$ et $\widetilde M\buildrel \gamma \over \longrightarrow \widetilde M$ la transformation affine $\gamma (z,t)=(\lambda z,\lambda t)$ ;
elle engendre une action analytique de $\Gamma=\langle \gamma \rangle =\{ \gamma^k:k\in {\Bbb Z}\} \simeq {\Bbb Z}$ sur $\widetilde M$. Cette action est libre et propre ; le quotient
$M=\widetilde M/\Gamma $
est donc une variété. Celle-ci est compacte, et plus précisément analytiquement isomorphe à la
variété de Hopf réelle ${\Bbb S}^n\times {\Bbb S}^1$. Les feuilletages $\widetilde{\cal F}$ et
$\widetilde{\cal V}$ sont invariants par $\gamma $ et induisent respectivement des feuilletages ${\cal F}$ et ${\cal V}$ sur $M$.
\smallskip
-- Les feuilles  $F_+$ et $F_-$ de ${\cal F}$ provenant respectivement  de
$\widetilde F_+$ et $\widetilde F_-$ sont des cercles (courbes fermées simples). Les autres sont des courbes non fermées ; pour chacune
d'elle, un bout spirale autour de $L_+$ et l'autre autour de $L_-$. Le feuilletage ${\cal F}$ est connu comme étant le {\sl flot affine de Reeb} sur
${\Bbb S}^n\times {\Bbb S}^1$.

-- La feuille $V$ de ${\cal V}$ provenant de $\widetilde V$ est une variété analytique isomorphe à la variété de Hopf
${\Bbb S}^{n-1}\times {\Bbb S}^1$. Les autres feuilles sont difféomorphes à ${\Bbb E}$. Le
feuilletage ${\cal V}$ est connu comme étant le {\sl feuilletage affine de Reeb} sur
${\Bbb S}^n\times {\Bbb S}^1$.
\smallskip
On peut mettre sur la variété ${\Bbb S}^n\times {\Bbb S}^1$ une métrique riemannienne $\Theta $ de telle sorte que les feuilletages ${\cal F}$ et
${\cal V}$ soient orthogonaux. Il suffit d'en définir une sur $\widetilde M$ invariante par $\gamma $, par exemple :
$$\Theta = {1\over {\vert \vert z\vert \vert^2+t^2}}\left( dz_1\otimes  dz_1+\cdots +dz_n\otimes dz_n+dt\otimes dt\right) \leqno{(8)}$$
où $\vert \vert z\vert \vert  =\sqrt{z_1^2+\cdots +z_n^2}$ est la norme euclidienne de $z$
dans ${\Bbb E}$.  Pour cette métrique, les deux feuilletages ${\cal F}$ et ${\cal V}$ sont conformes.
\medskip
{\zazou 3. Cohomologie feuilletée de $(M,{\cal F})$}
\medskip
\noindent Dans cette section on calcule explicitement la cohomologie feuilletée du feuilletage ${\cal F}$
sur la variété $M={\Bbb S}^n\times {\Bbb S}^1$.
Plus précisément, on démontre  le :
\smallskip
\noindent {\bf 3.1. Théorème.} {\it La cohomologie feuilletée $H_{\cal F}^\ell (M)$ du flot affine
de Reeb ${\cal F}$ sur la variété de Hopf $M={\Bbb S}^n\times {\Bbb S}^1$ est ${\Bbb C}$ pour $\ell =0$
et est isomorphe (en tant qu'espace de Fréchet) à  $C^\infty ({\Bbb S}^{n-1}\times {\Bbb S}^1)$ pour $\ell =1$. (Et bien sûr, $H_{\cal F}^\ell (M)=0$ pour $\ell \geq 2$.)}
\smallskip
 Pour $\ell \geq 2$, on a
 $H_{\cal F}^\ell (M)=0$ pour une raison évidente de degré. D'autre part, il est
 facile de voir que $H_{\cal F}^0(M)={\Bbb C}$. Cela vient du fait que toute feuille non compacte de ${\cal F}$ spirale positivement sur $F_+$ et
négativement sur $F_-$ ; donc une fonction basique (continue) est nécessairement constante.
\smallskip
Toute la suite de ce texte sera consacrée à  la démonstration de
l'isomorphisme
$H_{\cal F}^1(M)\simeq C^\infty ({\Bbb S}^{n-1}\times {\Bbb S}^1)$. Elle se fera en trois étapes.

\smallskip
\noindent {\bf 3.2.} Dans l'hyperplan affine de ${\Bbb E}\times {\Bbb R}$ d'équation $t=1$, on note $B_+$ la boule fermée de
centre $(0,1)$ et de rayon un nombre $\rho >0$ et $C_+$ le cône fermé de base $B_+$. De même, dans l'hyperplan affine d'équation $t=-1$, on considère
la boule $B_-$ de
centre $(0,-1)$ et de rayon $\rho >0$ et $C_-$ le cône fermé de base $B_-$. On pose :
$$U_+={\Bbb E}\times {\Bbb R}\setminus C_-,\hskip1cm U_-={\Bbb E}\times {\Bbb R}\setminus C_+\hskip1cm \hbox{et}\hskip1cm U=U_+\cap U_-.\leqno{(9)}$$
Alors ${\cal U}=\{ U_+,U_-\} $ est un recouvrement ouvert de $\widetilde M$. Les feuilletages induits par
$\widetilde{\cal F}$ sur $U_+$ et $U_-$ sont isomorphes au feuilletage sur ${\Bbb E}\times {\Bbb R}$ défini par
la première projection $(z,t)\longrightarrow z$. Ils sont donc {\sl intégrablement homotopes} ({\it cf.} [Ek])
au feuilletage par points sur l'espace euclidien ${\Bbb E}$. De la même manière, le feuilletage induit
par
$\widetilde{\cal F}$ sur $U$ est isomorphe au feuilletage sur ${\Bbb E}^\ast \times {\Bbb R}$ défini par
la première projection $(z,t)\longrightarrow z$. Il est donc {\sl intégrablement homotope}
au feuilletage par points sur ${\Bbb E}^\ast $. Leurs cohomologies feuilletées respectives sont donc triviales ; plus précisément on a :
$$H_{\widetilde{\cal F}}^\ell (U_+)=H_{\widetilde{\cal F}}^\ell (U_-)=\cases{C^\infty ({\Bbb E})&\hbox{si $\ell =0$}\cr
0&\hbox{si $\ell \geq 1$}}\leqno{(10)}$$
et
$$H_{\widetilde{\cal F}}^\ell (U)=\cases{C^\infty ({\Bbb E}^\ast )&\hbox{si $\ell =0$}\cr
0&\hbox{si $\ell \geq 1$}}\leqno{(11)}$$
Voici par exemple les ouverts $U_+$ et $U$ ($U_-$ est obtenu simplement en appliquant à $U_+$ la  réflexion par rapport à l'hyperplan d'équation $t=0$).
\medskip
\centerline{\hskip0.7cm\epsfxsize=8.5cm
\epsfbox{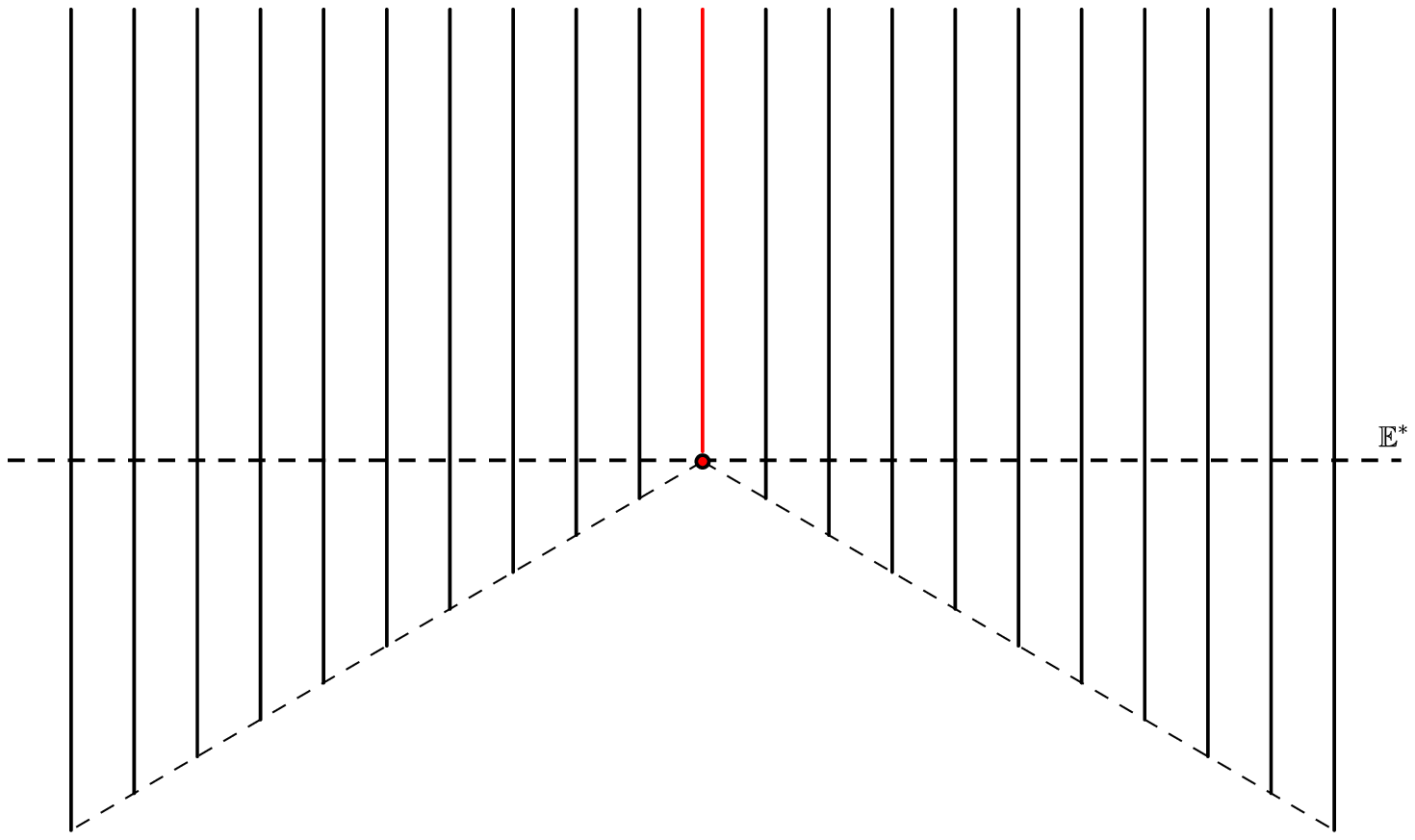}}
\smallskip
\centerline{L'ouvert $U_+$}

\bigskip
\centerline{\hskip0.7cm\epsfxsize=8.5cm
\epsfbox{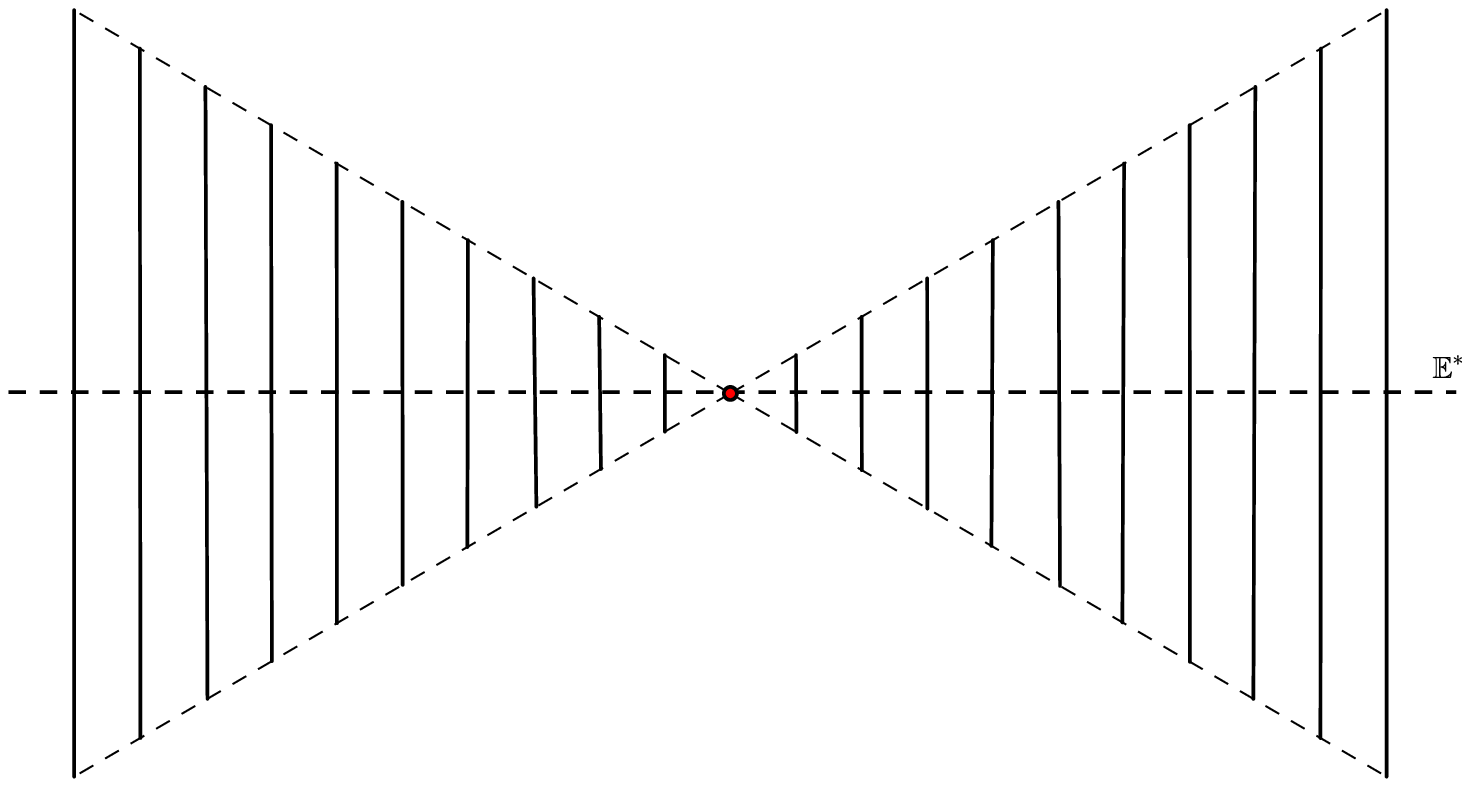}}
\smallskip
\centerline{L'ouvert $U$}
\bigskip
Le recouvrement ouvert ${\cal U}=\{ U_+,U_-\} $ est donc acyclique pour le faisceau ${\cal O}_{\widetilde{\cal F}}$. Par suite :
$$H_{\widetilde{\cal F}}^\ast (\widetilde M)=H^\ast (\widetilde M,{\cal O}_{\widetilde{\cal F}})=
H^\ast ({\cal U},{\cal O}_{\widetilde{\cal F}}).\leqno{(12)}$$
\smallskip
\noindent {\bf 3.3.} On  peut alors utiliser le recouvrement ouvert  ${\cal O}_{\widetilde{\cal F}}$-acyclique ${\cal U}$
pour calculer la cohomologie feuilletée de  $(\widetilde M,\widetilde{\cal F})$. La suite
 exacte longue de Mayer-Vietoris associée s'écrit :
$$0\longrightarrow H_{\widetilde{\cal F}}^0(\widetilde M)\buildrel r\over
\longrightarrow H_{\widetilde{\cal F}}^0(U_+)\oplus H_{\widetilde{\cal F}}^0(U_-)\buildrel j\over
\longrightarrow H_{\widetilde{\cal F}}^0(U)\buildrel \delta \over \longrightarrow H_{\widetilde{\cal F}}^1(\widetilde M)
\longrightarrow 0\leqno{(13)}$$
où $r(f)=\left( f_{\vert U_+},f_{\vert U_-}\right) $ (qu'on notera $(f,f)$),
$j(g,h)=g_{\vert U}-h_{\vert U}$ (qu'on notera $g-h$) et $\delta $ est l'homomorphisme de connexion. Cette suite exacte
s'écrit plus précisément :
$$0\longrightarrow C^\infty ({\Bbb E})\buildrel r\over \longrightarrow C^\infty ({\Bbb E})
\oplus C^\infty ({\Bbb E})\buildrel j\over \longrightarrow C^\infty ({\Bbb E}^\ast )\buildrel \delta \over \longrightarrow H_{\widetilde{\cal F}}^1(\widetilde M)
\longrightarrow 0\leqno{(14)}$$
toujours avec $r(f)=(f,f)$ et $j(g,h)=g-h$. L'image de $r$ est un sous-espace fermé de $C^\infty ({\Bbb E})\oplus C^\infty ({\Bbb E})$ isomorphe à
$C^\infty ({\Bbb E})$ ; il a pour supplémentaire le facteur $(C^\infty ({\Bbb E}),0)$ puisque tout $(g,h)$ s'écrit de façon unique $(g,h)=(h,h)+(g-h,0)$. En plus,
ce facteur est l'image de $j$ dans $C^\infty ({\Bbb E}^\ast )$. On a donc un isomorphisme canonique :
$$H_{\widetilde{\cal F}}^1(\widetilde M)\simeq C^\infty ({\Bbb E}^\ast )/C^\infty ({\Bbb E}).\leqno{(15)}$$
\smallskip
\noindent Pour simplifier, on notera dorénavant $W$ l'espace $H_{\widetilde{\cal F}}^1(\widetilde M)\simeq C^\infty ({\Bbb E}^\ast )/C^\infty ({\Bbb E}).$
Muni de la topologie quotient, $W$ hérite d'une structure d'EVT non séparé puisque $C^\infty ({\Bbb E})$ est un sous-espace non fermé de $C^\infty ({\Bbb E}^\ast )$
({\it cf.} Appendice section 5).
\smallskip
\noindent {\bf 3.4.} Comme l'ouvert ${\Bbb E}^\ast $  est invariant par l'homothétie $\gamma (z)=\lambda z$, $\gamma $ agit sur $C^\infty ({\Bbb E}^\ast )$
en préservant $C^\infty ({\Bbb E})$ donc aussi sur le quotient $W=C^\infty ({\Bbb E}^\ast )/C^\infty ({\Bbb E})$ ; par suite $W$ est un $\Gamma $-module. On a
donc une suite exacte de $\Gamma $-modules :
$$0\longrightarrow C^\infty ({\Bbb E})\buildrel j \over \hookrightarrow C^\infty ({\Bbb E}^\ast )\buildrel \tau \over \longrightarrow W\longrightarrow 0\leqno{(16)}$$
\medskip
\noindent où $j$ est l'application dans la suite (14) et $\tau $ la projection canonique de $C^\infty ({\Bbb E}^\ast )$ sur le quotient $W$. En prenant la cohomologie de $\Gamma $
à valeurs dans chacun de ces $\Gamma $-modules, on obtient  une suite exacte
longue :
\medskip
\centerline{$0\longrightarrow H^0(\Gamma ,C^\infty ({\Bbb E}))\buildrel j^\ast  \over \longrightarrow H^0(\Gamma ,C^\infty ({\Bbb E}^\ast ))\buildrel \tau^\ast  \over \longrightarrow
H^0(\Gamma ,W)\buildrel \delta  \over \longrightarrow \cdots $}
$$\cdots \buildrel \delta  \over \longrightarrow H^1(\Gamma ,C^\infty ({\Bbb E}))
\buildrel j^\ast  \over \longrightarrow H^1(\Gamma ,C^\infty ({\Bbb E}^\ast ))\buildrel \tau^\ast  \over \longrightarrow
H^1(\Gamma ,W)\longrightarrow 0\leqno{(17)}$$
où $\delta $ est l'homomorphisme de connexion habituel. Reste à déterminer les espaces en question, nous aurons bien besoin de certains d'entre eux pour la suite.

\smallskip
$\bullet $ L'espace $H^0(\Gamma ,C^\infty ({\Bbb E}))$ est celui des fonctions de classe  $C^\infty $ sur ${\Bbb E}$ invariantes par $\gamma $, donc constantes ; par suite
$H^0(\Gamma ,C^\infty ({\Bbb E}))={\Bbb C}$.

\smallskip
$\bullet $  L'espace $H^0(\Gamma ,C^\infty ({\Bbb E}^\ast ))$ est celui des fonctions $C^\infty $ sur ${\Bbb E}^\ast $ invariantes par $\gamma $, donc
les fonctions sur $V={\Bbb E}^\ast /\langle \gamma \rangle ={\Bbb S}^{n-1}\times {\Bbb S}^1$ {\it i.e.} :
$$H^0(\Gamma ,C^\infty ({\Bbb E}^\ast ))=C^\infty ({\Bbb S}^{n-1}\times {\Bbb S}^1).$$

\smallskip
$\bullet $  Calculons l'espace $H^1(\Gamma ,C^\infty ({\Bbb E}))$. C'est le quotient de $C^\infty ({\Bbb E})$ par le sous-espace engendré par les éléments
de la forme $f-f\circ \gamma $. Nous avons donc à résoudre l'équation cohomologique $f(z)-f(\lambda z)=g(z)$ pour $g\in C^\infty ({\Bbb E})$ donnée.
Une condition nécessaire est $g(z)=0$. Si on la suppose satisfaite, une solution formelle est donnée par la série :
$$f(z)=\displaystyle \sum_{k=0}^\infty g(\lambda^kz).$$
Cette série converge pour $z=0$, puisque tous ses termes sont nuls. D'autre part, on montre facilement que, pour tout $R>0$, toutes ses séries dérivées sont équivalentes, sur
la boule fermée de rayon $R$ et centrée à l'origine, à des
séries géométriques de raison une puissance de $\lambda $, donc elles convergent. Finalement, cette série converge pour la topologie $C^\infty $. Ce qui montre que $f$ est
$C^\infty $ et est une solution de l'équation cohomologique $f-f\circ \gamma =g$. On en déduit que $H^1(\Gamma ,C^\infty ({\Bbb E}))$ est de dimension $1$ sur ${\Bbb C}$ engendré
par la fonction constante égale à $1$.
\smallskip
$\bullet $ Pour  calculer $H^1(\Gamma ,C^\infty ({\Bbb E}^\ast ))$, on utilise le $\Gamma $-revêtement $\sigma :{\Bbb E}^\ast \longrightarrow {\Bbb S}^{n-1}\times {\Bbb S}^1$
où ${\Bbb S}^{n-1}\times {\Bbb S}^1$ est le quotient de ${\Bbb E}^\ast $ (vu comme
la feuille $\widetilde V={\Bbb E}^\ast \times \{ 0\} $ du feuilletage ${\cal V}$ décrit dans la sous-section 2.1)
par l'action du groupe $\Gamma =\langle \gamma \rangle $. Notons  ${\cal C}^\ast $ et ${\cal C}$ les faisceaux des germes de fonctions $C^\infty $ respectivement
sur ${\Bbb E}^\ast $ et $V={\Bbb S}^{n-1}\times {\Bbb S}^1$. Ces faisceaux sont fins et donc :
$$H^\ell ({\Bbb E}^\ast ,{\cal C}^\ast )=\cases{C^\infty ({\Bbb E}^\ast ) & \hbox{ si $\ell =0$}\cr 0 &\hbox{ si $\ell \geq 0$}}
\hskip0.5cm \hbox{et}\hskip0.5cm
H^\ell (V,{\cal C})=\cases{C^\infty (V)& \hbox{ si $\ell =0$}\cr 0 & \hbox{ si $\ell \geq 0$}.}\leqno{(18)}$$
\medskip
\noindent Au revêtement $\sigma $ est associée une suite
spectrale convergeant vers $H^{p+q}(V,{\cal C})$ et de terme $D_2^{pq}=H^p(\Gamma ,H^q({\Bbb E}^\ast , {\cal C}^\ast ))$. Comme :
$$0=H^1(V,{\cal C})=H^1(\Gamma ,C^\infty ({\Bbb E}^\ast ))\oplus H^0(\Gamma ,H^1({\Bbb E}^\ast ,{\cal C }^\ast)),$$
on a $H^1(\Gamma ,C^\infty ({\Bbb E}^\ast ))=0$.

\smallskip
$\bullet $ Avec ce qu'on vient de calculer, la suite exacte longue de cohomologie (17) peut donc s'écrire  :
$0\longrightarrow {\Bbb C}\hookrightarrow C^\infty ({\Bbb S}^{n-1}\times {\Bbb S}^1)\longrightarrow H^0(\Gamma ,W)\buildrel \delta \over \longrightarrow {\Bbb C}
\longrightarrow 0.$
Elle donne la suite exacte (plus courte) :
$$0\longrightarrow C^\infty ({\Bbb S}^{n-1}\times {\Bbb S}^1)/{\Bbb C}\longrightarrow H^0(\Gamma ,W)\buildrel \delta \over \longrightarrow {\Bbb C}
\longrightarrow 0\leqno{(19)}$$ qui montre que :
$$H^0(\Gamma ,W)=\left( C^\infty ({\Bbb S}^{n-1}\times {\Bbb S}^1)/{\Bbb C}\right) \oplus {\Bbb C}\simeq C^\infty ({\Bbb S}^{n-1}\times {\Bbb S}^1).\leqno{(20)}$$

$\bullet $ Pour finir, la suite spectrale associée au revêtement feuilleté $(\widetilde M,\widetilde{\cal F})\longrightarrow (M,{\cal F})$ donne :
$H_{\cal F}^1(M)=H^0(\Gamma ,H_{\widetilde{\cal F}}^1(\widetilde M))\oplus H^1(\Gamma ,H_{\widetilde{\cal F}}^0(\widetilde M)).$
Mais, comme on vient de le voir dans tous les calculs qui précèdent, on a :
$$H^1(\Gamma ,H_{\widetilde{\cal F}}^0(\widetilde M))=H^1(\Gamma ,C^\infty (\Bbb E^\ast ))=0$$
et :
$$H^0(\Gamma ,H_{\widetilde{\cal F}}^1(\widetilde M))=H^0(\Gamma ,W)\simeq C^\infty ({\Bbb S}^{n-1}\times {\Bbb S}^1).$$
Finalement on a le résultat annoncé $H_{\cal F}^1(M) \simeq C^\infty ({\Bbb S}^{n-1}\times {\Bbb S}^1)$. Ce qui termine la démonstration du théorème 3.1. \hfill $\diamondsuit  $
\smallskip
\noindent {\bf 3.5. Remarque}

\smallskip
Le théorème 3.1. dit que l'espace de cohomologie feuilletée $H_{\cal F}^1({\Bbb S}^n\times {\Bbb S}^1)$ est paramétré par les fonctions de classe $C^\infty $
sur la feuille compacte $V={\Bbb S}^{n-1}\times {\Bbb S}^1$ du feuilletage orthogonal ${\cal V}$. Cette feuille est presque une section de ${\cal F}$ : elle coupe une
et seule fois chacune des feuilles de ${\cal F}$ sauf $F_+$ et $F_-$. Si on prive $M$ de $F_+$ et $F_-$, on obtient une variété $M_0$ munie du feuilletage induit
${\cal F}_0$. Le revêtement universel $\widetilde M_0$ de $M_0$ est le produit ${\Bbb E}^\ast \times {\Bbb R}$ et le feuilletage $\widetilde {\cal F}_0$ relevé a pour feuilles
les facteurs $\{ z\} \times {\Bbb R}$ avec $z\in {\Bbb E}^\ast $.
Si on munit $\widetilde V_0={\Bbb E}^\ast \times \{ 0\} $ du feuilletage par points (qu'on notera ${\cal P}_0$), on a une rétraction feuilletée
$H((z,t),s)=(z,st)$ (avec $s\in [0,1]$) de $(\widetilde M_0,\widetilde {\cal F}_0)$ sur $(\widetilde V_0,{\cal P}_0)$. Elle est équivariante par rapport à l'action
de $\gamma $ mais ne préserve pas chaque feuille individuellement ; elle ne saurait donc être une rétraction intégrable. Toutefois, elle le devient quand on prend les
quotients respectifs sous l'action de $\Gamma $. Ainsi, le feuilletage ${\cal F}_0$ sur $M_0$ est intégrablement homotope au feuilletage par points sur
$V={\Bbb S}^{n-1}\times {\Bbb S}^1$. Donc $H_{{\cal F}_0}^1(M_0)=0$.

\medskip
{\zazou 4. Distributions invariantes}
\medskip
\noindent {\bf 4.1.} Si ${\cal F}$ est un flot défini par un champ non singulier $X$ sur une variété $M$, sa cohomologie feuillet\'ee peut se voir comme suit.
Soit $\nu $ un sous-fibr\'e suppl\'ementaire \`a ${\cal F}$
dans $TM$. Soit $\chi $ la $1$-forme diff\'erentielle telle que
$\chi (X)=1$ et $\chi_{\vert \nu }=0$. Il est facile de voir que,
pour tout $\ell \in {\Bbb N}$, on a :
$$\Omega_{\cal F}^\ell (M)=\cases{C^\infty (M)&\hbox{si $\ell =0$}\cr
C^\infty (M)\otimes \chi &\hbox{si $\ell =1$}\cr 0&\hbox{si $\ell \geq
2$}}\leqno{(21)}$$ et que le   complexe feuillet\'e  se r\'eduit \`a :
$0\longrightarrow \Omega_{\cal F}^0(M)\buildrel {d_X}\over
\longrightarrow \Omega_{\cal F}^1(M)\longrightarrow 0$ o\`u $d_X$
est l'op\'erateur   d\'efini par $d_Xf=(X\cdot f)\otimes \chi $. Son
conoyau $\Omega_{\cal F}^1(M)/\hbox{Im}d_X$ est exactement le
premier espace de {\it cohomologie feuillet\'ee} $H_{\cal F}^1(M)$
de ${\cal F}$. Il ne d\'epend pas du champ qui le d\'efinit : on
v\'erife ais\'ement, en exhibant explicitement un isomorphisme de
complexes feuillet\'es, qu'on obtient la m\^eme cohomologie si
on remplace le champ $X$ par un champ $Z=hX$ avec $h$ fonction
partout non nulle.

Le calcul de l'espace $H_{\cal F}^1(M)$ revient  \`a la
r\'esolution de l'\'equation cohomologique continue pour le champ
$X$ ({\it cf.} [DE]) :
$$\hbox{\it \'Etant donnée $g\in C^\infty (M)$, existe-t-il
 $f\in C^\infty (M)$ telle que $X\cdot f=g$ ?}\leqno{(22)}$$
L'espace $H_{\cal F}^1(M)$
contient exactement les obstructions à la résolution de cette équation, d'où l'intérêt de son calcul.
\smallskip
\noindent {\bf 4.2.} Supposons $M$ compacte. Une {\it distribution} sur $M$ est une forme lin\'eaire continue :
$$\varphi \in C^\infty (M)\buildrel T\over \longrightarrow
\langle T,\varphi \rangle \in {\Bbb C}$$   {\it i.e.} un \'el\'ement
du dual topologique  ${\cal D}'(M)$ de $C^\infty (M)$. L'espace
vectoriel ${\cal D}'(M)$ sera muni de la {\it topologie faible} {\it
i.e.} la topologie la moins fine qui rend continues toutes les
\'evaluations lin\'eaires $e_\varphi :T\in {\cal D}'(M)\longmapsto
\langle T,\varphi \rangle \in {\Bbb C}$. Toute  forme volume $\mu $
sur $M$ d\'efinit une injection $f\in
C^\infty (M)\longmapsto T_f\in {\cal D}'(M)$ donn\'ee par $\langle
T_f,\varphi \rangle =\int_M(f\varphi )\mu .$ Une distribution de
ce type est dite {\it r\'eguli\`ere}.
\smallskip

Soit $M\buildrel \gamma \over \longrightarrow M$ un
diff\'eomorphisme. Une distribution $T$ sur $M$ est dite
{\it invariante} par $\gamma $ (ou simplement $\gamma $-{\it
invariante}) si elle v\'erifie $\langle T,\varphi \circ \gamma
\rangle =\langle T,\varphi \rangle $ pour toute fonction $\varphi
\in C^\infty (M)$.   On dira que $T$ est invariante par un groupe
$\Gamma $ de diff\'eomorphismes de $M$ (ou $\Gamma $-invariante) si
elle est invariante par chacun de ses \'el\'ements.
\smallskip
\noindent {\bf 4.3.}
Soit $X$ un champ de vecteurs  sur $M$. C'est aussi un opérateur différentiel du premier ordre (on l'a déjà considéré ainsi)
 $X: C^\infty (M)\longrightarrow C^\infty (M)$ défini par l'égalité $(X\cdot
f)(x)=d_xf(X_x)$.
Il admet une
extension naturelle aux distributions :
$$X:T\in {\cal
D}'(M)\longrightarrow X\cdot T\in {\cal D}'(M)$$ avec $\langle X\cdot
T,\varphi \rangle = -\langle T,X\cdot \varphi \rangle $. (On peut
donc s'int\'eresser aussi \`a la r\'esolution de l'\'equation
cohomologique continue au niveau des distributions
$X\cdot T=S$.)
\smallskip
Une distribution $T$ est dite {\it invariante} par $X$ ou $X$-{\it
invariante} si elle v\'erifie $X\cdot T=0$ {\it i.e.} elle est nulle
sur l'image de  $X:C^\infty (M)\longrightarrow C^\infty (M)$ qui est
l'espace des {\it divergences} de $X$. Une condition n\'ecessaire
(et non suffisante en g\'en\'eral) pour que l'\'equation (22) admette
une solution $f$ est $\langle T,g\rangle =0$ pour toute
distribution $T$ invariante par $X$. L'espace ${\cal D}_X'(M)$ des distributions invariantes par $X$ s'identifie donc naturellement
au dual topologique de l'espace $H_{\cal F}^1(M)$.
\smallskip
Revenons à notre flot ${\cal F}$ sur la variété de Hopf $M={\Bbb S}^n\times {\Bbb S}^1$. Il peut être défini
par n'importe quel champ non singulier dont les feuilles sont les courbes intégrales. Un tel champ s'écrit en coordonnées sur $\widetilde M$ sous la forme
$\widetilde X=a(z,t){\partial \over {\partial t}}$ où $a$ est une fonction $C^\infty $ partout non nulle et vérifiant la condition d'invariance $a(\lambda z,\lambda t)=\lambda a(z,t)$,
par exemple $a(z,t)=\sqrt{\vert \vert z\vert \vert^2+t^2}$. Du théorème 3.1 on tire donc le :
\smallskip
\noindent {\bf 4.4. Corollaire.} {\it Pour un champ $X$ définissant le feuilletage ${\cal F}$, l'espace ${\cal D}_X'(M)$ des
distributions $X$-invariantes est naturellement isomorphe au dual topologique de $H_{\cal F}^1({\Bbb S}^n\times {\Bbb S}^1)=C^\infty ({\Bbb S}^{n-1}\times {\Bbb S}^1)$,
donc à l'espace des distributions ${\cal D}'({\Bbb S}^{n-1}\times {\Bbb S}^1)$ sur la variété de Hopf  ${\Bbb S}^{n-1}\times {\Bbb S}^1$.}

\medskip
{\zazou 5. Appendice}
\medskip

\noindent Il est bien connu que le premier groupe de cohomologie feuilletée d'un flot linéaire sur le tore ${\Bbb T}^2$ à pente un irrationnel de Liouville
est un espace vectoriel topologique non séparé. Et souvent, quand la question se pose à cet effet, curieusement
tout le monde cite cet exemple ! Mais ça peut se produire même lorsque les feuilles sont fermées
et que le feuilletage est presque une fibration, c'est le cas du feuilletage $(\widetilde M,\widetilde{\cal F})$ : on a  vu
que $H_{\widetilde{\cal F}}^1(\widetilde M)$ s'identifie au quotient $C^\infty ({\Bbb E}^\ast )/C^\infty ({\Bbb E})$ dont nous avons affirmé qu'il n'est pas séparé. L'objet
de cet appendice est de justifier cela en montrant que l'espace $C^\infty ({\Bbb E})$ n'est pas fermé dans $C^\infty ({\Bbb E}^\ast )$.
\smallskip
\noindent {\bf 5.1. L'espace de Fréchet $C^\infty ({\Bbb E}^\ast )$}
\smallskip
Soient $\lambda $ un nombre réel tel que $\lambda >1$. Pour tout entier $k\in {\Bbb N}^\ast $, on note $C_k$ la couronne de
${\Bbb E}^\ast $ donnée par $C_k=\{ z\in {\Bbb E}:\lambda^{-k}\leq \vert \vert z\vert \vert \leq \lambda^k \} $.
Alors la famille $\{  C_k\} $ est une suite croissante de compacts telle que
$\displaystyle \bigcup_{k=1}^\infty C_k={\Bbb E}^\ast .$
\smallskip
Pour un multi-indice $s=(s_1,\cdots ,s_n)\in {\Bbb N}^n$, on note $\vert s\vert =s_1+\cdots +s_n$ sa longueur et $D^s$ l'opérateur différentiel ${{\partial^{\vert s\vert }}\over \partial^{s_1}z_1\cdots \partial^{s_n}z_n}$. Pour $k,r\in {\Bbb N}$ et $f\in C^\infty ({\Bbb E}^\ast )$, on pose :
$$\rho_{k,r}(f)=\sum_{\vert s\vert \leq r } \sup_{z\in C_k}\left\vert D^sf(z)\right\vert .\leqno{(23)}$$
On obtient ainsi une famille dénombrable filtrante de semi-normes sur l'espace $ C^\infty ({\Bbb E}^\ast )$. Elle  y définit une topologie ${\cal T}$ qui en fait un espace de Fréchet.
Une suite de fonctions $f_p$ dans $ C^\infty ({\Bbb E}^\ast )$ converge vers $f\in  C^\infty ({\Bbb E}^\ast )$ au sens de cette topologie si, pour tous $k,r\in {\Bbb N}$,
la suite numérique $\rho_{k,r}(f_p-f)$ (indexée par $p$) tend vers $0$.

\smallskip
\noindent {\bf 5.2. Le sous-espace $C^\infty ({\Bbb E})$ n'est pas fermé dans $C^\infty ({\Bbb E}^\ast )$}
\smallskip
Pour le montrer, nous allons exhiber explicitement une suite de fonctions dans $C^\infty ({\Bbb E})$ qui converge au sens de la topologie ${\cal T}$ vers
une fonction $f$ dans $C^\infty ({\Bbb E}^\ast )\setminus C^\infty ({\Bbb E})$.
\smallskip
Soit $\phi :]0,+\infty [\longrightarrow {\Bbb R}$ une fonction de classe $C^\infty $ ne se prolongeant pas en fonction $C^\infty $ sur $[0,+\infty [$.
(On peut penser par exemple à $\phi (t)=\sqrt{t}$, $\phi (t)={1\over t}$ ou...) On pose :
$$f(z) = \phi (\vert \vert z\vert \vert^2)\hskip0.3cm \hbox{et}\hskip0.3cm f_p(z)=\phi \left( \vert \vert z\vert \vert^2 +{1\over p}\right)\hskip0.3cm\hbox{pour $p\in {\Bbb N}^\ast $}.\leqno{(25)}$$
Il est évident que, pour tout $p\geq 1$, la fonction $f_p$ est dans $C^\infty ({\Bbb E})$ ; par contre $f$ est dans $C^\infty ({\Bbb E}^\ast )$ mais pas dans $C^\infty ({\Bbb E})$.
\smallskip
On va indiquer rapidement comment montrer que la suite $(f_p)$ converge vers $f$ pour la topologie ${\cal T}$. Un calcul facile mais lourd montre que, pour tout $s=(s_1,\cdots ,s_n)\in {\Bbb N}^n$, on a :
$$D^sf(z)=\sum_{\ell =0}^{\vert s\vert }a_\ell (z)\phi^{(\ell )}(\vert \vert z\vert \vert^2)\hskip0.5cm \hbox{et}\hskip0.5cm D^sf_p(z)=\sum_{\ell =0}^{\vert s\vert }a_\ell (z)\phi^{(\ell )}\left( \vert \vert z\vert \vert^2+{1\over p}\right) $$
où $a_\ell $ est une fonction polynomiale en $z_1,\cdots ,z_n$ et $\phi^{(\ell )}$ désigne la dérivée d'ordre $\ell $ de la fonction $\phi $ (sur $]0,+\infty [$ bien sûr).
On se met sur l'un des compacts $C_k$.
Là, encore une fois, un calcul
facile mais un peu long, usant du théorème des accroissements finis appliqué à chacune des fonctions $\phi^{(\ell )}$ sur l'intervalle $\left\lbrack \vert \vert z\vert \vert^2,\vert \vert z\vert \vert^2+{1\over p}\right\rbrack $ donne :
$$\sup_{z\in C_k}\left\vert D^sf_p(z)-D^sf(z)\right\vert  \leq \left\{ \sum_{\ell =0}^{\vert s\vert }\left( \sup_{z\in C_k}\vert a_\ell (z)\vert \right) \left( \sup_{z\in C_k}\vert \phi^{(\ell +1)}(z)\vert \right) \right\} {1\over p}.$$
D'où :
$$\rho_{k,r}(f_p-f)\leq \left\{ \sum_{\vert s\vert \leq r}\sum_{\ell =0}^{\vert s\vert }\left( \sup_{z\in C_k}\vert a_\ell (z)\vert \right) \left( \sup_{z\in C_k}\vert \phi^{(\ell +1)}(z)\vert \right) \right\} {1\over p}\leq {C\over p}\leqno{(26)}$$
où $C$ est une constante strictement positive qui ne dépend que de $\phi $, $k$ et $r$. Ceci
montre que $f_p$ converge vers $f$ pour toute semi-norme $\rho_{k,r}$. Autrement dit, la suite $(f_p)$ converge vers $f$ pour la topologie ${\cal T}$. \hfill $\diamondsuit  $

\vfill\break

\centerline{\tfp  R\'ef\'erences}
\vskip0.3cm
\item {[DE]} {\sc  Deghan-Nezhad, A. \& El Kacimi Alaoui, A.} {\it \'Equations cohomologiques de flots riemanniens et de difféomorphismes  d'Anosov}.
Journal of the Math. Society of Japan, Vol. 59 N° 4. (2007), 1105-1134.
\smallskip
\item{[Ek]} {\sc  El Kacimi Alaoui, A.} {\it  Sur la cohomologie feuillet\'ee.}
Compositio Mathematica 49, (1983), 195-215.

\smallskip

\noindent [Go] {\sc  Godbillon, C.} {\it Feuilletages - \'Etudes géométriques.} Birkh\"auser, (1991).

\smallskip
\item{[MM]} {\sc  Matsumoto, S. \& Mitsumatsu, Y.} {\it Leafwise cohomology and rigidity of certain Lie group actions.} Ergod. Th. \& Dynam. Sys. 23, (2003), 1839-1866.
\smallskip

\item{[Su]} {\sc  Sullivan, D.} {\it  Cycles for the dynamical study of foliated manifolds
and complex mani\-folds.}  Invent. Math. 36, (1978), 225-255.

\smallskip

\item{[Va]} {\sc  Vaisman, I.} {\it  Cohomology and Differential Forms.}  M. Dekker, (1973).

\vskip1.5cm

\noindent {\pg Universit\'e Polytechnique Hauts-de-France

\noindent EA 4015 - LAMAV, FR CNRS 2956

\noindent F-59313 Valenciennes Cedex 9

\noindent FRANCE}

\smallskip

\noindent aziz.elkacimi@uphf.fr

\noindent http://perso.numericable.fr/azizelkacimi/

\medskip

\end